\documentclass[12pt]{article}

\usepackage{amsmath, amssymb,amsthm,latexsym}

\usepackage[dvips]{graphics}
\usepackage{epsfig, rotate}
\usepackage[english]{babel}
\usepackage{amsfonts}
\usepackage{array}
\usepackage{color}
\usepackage{float}

\newcommand\blfootnote[1]{%
  \begingroup
  \renewcommand\thefootnote{}\footnote{#1}%
  \addtocounter{footnote}{-1}%
  \endgroup
}

\newtheorem{theorem}{{Theorem}}

\newtheorem{lemma}{{Lemma}}

\numberwithin{theorem}{section}
\numberwithin{lemma}{section}
\numberwithin{corollary}{section}

\usepackage
{graphics}
\usepackage
{graphicx}

\newcommand{\be}{\begin{equation}}
\newcommand{\ee}{\end{equation}}
\newcommand{\beaa}{\begin{eqnarray*}}
\newcommand{\eeaa}{\end{eqnarray*}}
\newcommand{\bea}{\begin{eqnarray}}
\newcommand{\eea}{\end{eqnarray}}

\newcommand{\bei}{\begin{itemize}}
\newcommand{\eei}{\end{itemize}}

\newcommand{\td}{\overset{d}\to}
\newcommand{\tp}{\overset{p}\to}

\def\E{\mathbb{E}}

\def\td{\stackrel{d}{\rightarrow}}

\begin{document}

\thispagestyle{empty}

\noindent {\bf \Large Inference of high quantiles of a
heavy-tailed\\
distribution from block data}

\vspace{20pt}
\noindent {\bf Yongcheng Qi$^{1,\,2}$, Mengzi Xie$^1$,
Jingping Yang$^3$}

\vspace{20pt}
\date{}

{\small

\noindent $^1$Department of Mathematics and Statistics, University of Minnesota Duluth,\\
1117 University Drive, Duluth, MN 55812, USA.\\
$^2$ Corresponding author. 

\vspace{5pt}

\noindent $^3$Department of Financial Mathematics, School of
Mathematical Sciences, \\Peking University, Beijing 100871, PR
China.

\blfootnote{\small \hspace{-15pt}$^*$This is an Accepted Manuscript
version of the following article, accepted for publication in
\textit{Statistics} [https://doi.org/10.1080/02331888.2023.2228442].
It is deposited under the terms of the Creative Commons
Attribution-NonCommercial License
(http://creativecommons.org/licenses/by-nc/4.0/), which permits
non-commercial re-use, distribution, and reproduction in any medium,
provided the original work is properly cited.}

\vspace{20pt}

\noindent{\bf Abstract.} In this paper we consider the estimation
problem for high quantiles of a heavy-tailed distribution from block
data when only a few largest values are observed within blocks. We
propose estimators for high quantiles and prove that these
estimators are asymptotically normal. Furthermore, we employ
empirical likelihood method and adjusted empirical likelihood method
to constructing the confidence intervals of high quantiles. Through
a simulation study we also compare the performance of the normal
approximation method and the adjusted empirical likelihood methods
in terms of the coverage probability and length of the confidence
intervals.


\vspace{20pt}


\vspace{20pt} \noindent {\bf 2020 Mathematics Subject
Classification:}~   62G32 }


\newpage

\section{Introduction}\label{intro}

In the past four decades, extreme value theory has been developed
and used to model rare events in many disciplines so as to assess
the risks of those events. In classical settings, the distributions
under the general framework of extreme value theory include a very
large class of semi-parametric distributions, and estimators on the
tail index and high quantiles of these distributions can be found in
the literature, see, e.g., Hill~\cite{Hill1975},
Pickands~\cite{Pickands1975}, Dekkers and de Haan~\cite{DH1989},
Dekkers et al.~\cite{DEH1989}, de Haan and Rootz\'en~\cite{DR1993},
Ferreira et al.~\cite{FDP2003}, Gong et al.~\cite{GLPY2015}, and {
Allouche et al~\cite{AEG2022}, among others.}

In this paper, we are interested in heavy-tailed distributions which
are widely used to model data in fields like meteorology, hydrology,
climatology, environmental science, telecommunications, insurance
and finance. We refer the readers to Embrechts et al.~\cite{ELM1997}
for details. Based on a random sample with $n$ independent and
identically distributed (iid) data points from a heavy-tailed
distribution, Hill's estimator for the tail index and Weissman's
estimators for high quantiles, proposed by Hill~\cite{Hill1975} and
Weissman~\cite{Weissman1978}, respectively, are two popular ones in
the literature. These estimators are based on only a few of the
largest observations. For some recent methodologies for heavy-tailed
distributions, see, e.g., Peng and Qi~\cite{PengQi2017} and
Paulauskas and Vaiciulis~\cite{PV2013,PV2017}.


We consider the inference on heavy-tailed distributions from
incomplete data in this paper. In real world, some data are
naturally divided into several groups or blocks, and only a small
proportion of the largest observations within blocks are available
for analysis.  For example, for rainfall or snowfall, fire losses,
frequently, only a few of yearly largest observations are accessible
publicly. See Qi~\cite{Qi2010} for more examples.

For the block data, Paulauskas~\cite{Paulauska2003} introduced their
estimators for the tail index of a heavy-tailed distribution based
on the ratios of the first largest and second largest observations
within blocks. Devydov et al.~\cite{DPR2000} also used the similar
idea to estimate the index of stable random vectors. Later on,
Qi~\cite{Qi2010} proposed some Hill-type estimators for block data,
as given in \eqref{newest}, which have the smaller asymptotic
variance. Qi~\cite{Qi2010} also employed the empirical likelihood
methods to constructing confidence intervals for the tail index.
More recent development can be found in the literature; see, e.g.,
Paulauskas and Vaiciulis~\cite{PV2011, PV2012},
Vaiciulis~\cite{Vai2012, Vai2014},  Xiong and Peng~\cite{XP2020},
and Hu et al.~\cite{HPS2022}.

In this paper,  we propose estimators for high quantiles of the
heavy-tailed distribution from the block data and to construct
confidence intervals under the same settings as in Qi~\cite{Qi2010}.
We also employ empirical likelihood method and adjusted empirical
likelihood method to construct confidence intervals. It is worth
mentioning that under  the complete data setting for heavy-tailed
distributions, several papers have applied empirical likelihood
methods to construct confidence intervals for the tail index and
high quantiles, see, e.g.,  Lu and Peng~\cite{LP2002}, Peng and
Qi~\cite{PengQi2006a, PengQi2006b}.

The rest of the paper is organized as follows. In
section~\ref{main}, we introduce our estimators for high quantiles
and present their limiting distribution.  In section~\ref{elk}, we
employ the empirical likelihood method to construct confidence
intervals for the logarithm of high quantiles. In
section~\ref{simulation}, we conduct a simulation study to compare
the confidence intervals based on the normal approximation of our
estimators and the empirical likelihood method. Finally, we give all
the proofs in section~\ref{proofs}.

\section{Estimators of high quantiles}\label{main}

 A cumulative distribution function $F$ is a  heavy-tailed distribution function if it
satisfies the following condition
\begin{equation}\label{F}
1-F(x)=x^{-1/\gamma}L(x)~~ \mbox{ for } x>0,
 \end{equation}
where $\gamma>0$ is an unknown parameter and  $1/\gamma$ is called
the tail index of the distribution function $F$, and $L$ is a slowly
varying function at infinity satisfying
\[
\lim_{t\to\infty}\frac{L(tx)}{L(t)}=1,~~x>0
\]

Assume a random sample of size $n$ from $F$ is available.  A $p$-th
high quantile of distribution $F$, denoted as $x_p$, is defined as
the $(1-p)$-th quantile of $F$, i.e., $x_p=F^-(1-p)$, where $F^-$
denotes the generalized inverse of $F$ and $p=p_n\in (0,1)$
satisfying $\lim_{n\to\infty}p_n=0$ and $\lim_{n\to\infty}np_n=c\in
[0, \infty)$.

The inference of index $\gamma$ and high quantile $x_p$ has
attracted much attention in past fifty years. When a full sample
$X_1, \cdots, X_n$ is available, Hill's estimator for $\gamma$ and
Weissman's estimator for $x_p$ are well known in the literature; see
Hill~\cite{Hill1975} and Weissman~\cite{Weissman1978}. Some recent
developments on constructing confidence intervals of $\gamma$ and
$x_p$ based on normal approximations and empirical likelihood
methods can be found in Peng and Qi~\cite{PengQi2017}.



In this paper, we consider the case when the data is not fully
available. We describe our settings in this paper as follows.
Without loss of generality we can always assume $X_i\ge 1$ for $1\le
i\le n$, otherwise we can simply replace $X_i$ by $\max(X_i,1)$.
First, divide the sample $X_1, \cdots, X_n$ into $k_n$ blocks (or
groups), $V_1, \cdots, V_{k_n}$, and each block contains
$m=m_n=[n/k_n]$ observations, where $[x]$ denotes the integer part
of $x>0$. To be more specific, $V_i=\{X_{(i-1)m+1}, \cdots,
X_{im}\}$ for $1\le i\le k_n$. Let $X_{m,1}^{(i)}\ge \cdots \ge
X_{m,m}^{(i)}$ denote the order statistics of the $m$ observations
in the $i$-th block.

Let $r\ge 1$ be an integer.  Now we assume that the $r+1$ largest
random variables within each of the $k_n$ blocks are observed, that
is, only the observations $\{X_{m, j}^{(i)}: j=1, \cdots, r+1,
i=1,\cdots, k_n\}$ are available for inference.  {   From the data
within $i$-th block, Hill's estimator for $\gamma$ can be defined as
$\frac1r\sum^{r}_{j=1}(\log X_{m,j}^{(i)}-\log X_{m,r+1}^{(i)})$ for
$i=1, \cdots, k_n$. By using the average of all $k_n$ Hill's
estimators,} Qi~\cite{Qi2010} proposed the following estimator for
$\gamma$:
\begin{equation}\label{newest}
{\widehat{\gamma}}_n=\frac{1}{k_nr}\sum^{k_n}_{i=1}\sum^{r}_{j=1}(\log
X_{m,j}^{(i)}-\log X_{m,r+1}^{(i)}),
\end{equation}
where $k_n$ satisfies
\[
k_n\to\infty ~~\mbox{and }~~\frac{k_n}{n}\to 0 ~~\mbox{ as }
n\to\infty.
\]
In order to extend the current setting conveniently, we express the
above condition as
\[
k_n\to\infty ~~\mbox{and }~~m_n\to \infty ~~\mbox{ as } n\to\infty.
\]
Note that $k_n$ and $m_n$ are the number of blocks and the number of
observations within each block, respectively, and $n\sim k_nm_n$ is
approximately the total number of observations in all $k_n$ blocks.


To investigate the limiting distributions for the estimator, a
condition stronger than \eqref{F} is required. Throughout this
paper, we assume that
 there exists a measurable function $A(t)$ with $\lim_{t\to\infty}A(t)=0$ such that
 \begin{equation}\label{U}
\lim_{t\to\infty}\frac{U(tx)/U(t)-x^{\gamma}}{A(t)}=x^{\gamma}\frac{x^{\rho}-1}{\rho}
 \end{equation}
for all $x>0$, where $U(y)=F^{-}(1-\frac1y)$ is the inverse function of $\frac
1{1-F}$ and $\rho< 0$. This condition is more general than
the following condition
\begin{equation}\label{2terms}
1-F(x)=cx^{-1/\gamma}+dx^{-\beta}+o(x^{-\beta}) ~~\mbox{ as
}x\to\infty,
\end{equation}
where $0<\gamma^{-1}<\beta\le \infty$, which is used in
Paulauskas~\cite{Paulauska2003}. In fact, if (\ref{2terms}) holds,
then one can verify that (\ref{U}) holds with
$A(t)=-\gamma(\beta\gamma-1)dc^{-\beta\gamma}t^{1-\beta\gamma}$ and
$\rho=1-\beta\gamma$.

The asymptotic normality of ${\widehat{\gamma}}_n$ is as follows.

\begin{theorem}\label{thm1} (Qi~\cite{Qi2010}) Assume (\ref{U}) holds.  If
 \begin{equation}\label{allkn}
k_n\to\infty, ~~ m_n\to \infty \mbox{ and }k_n^{1/2}A(m_n)\to
\delta\in (-\infty,\infty) ~~\mbox { as } n\to\infty,
 \end{equation}
then
\[
(rk_n)^{1/2}({\widehat{\gamma}}_n-\gamma){\stackrel{d}{\to}}
N(\delta b_r, \gamma^2),
\]
where
\begin{equation}\label{br}
b_r=\frac{1}{r\rho}\left(\sum^r_{j=1}\frac{\Gamma(j-\rho)}{(j-1)!}-
\frac{\Gamma(r+1-\rho)}{(r-1)!}\right)
\end{equation}
and $\Gamma(x)=\int^{\infty}_0t^{x-1}e^{-t}dt$ is the Gamma
function.
\end{theorem}




\vspace{10pt}

In this paper, we propose the following estimator for $x_p$:
\[
\widehat{x}_p=\exp\Big(\frac{1}{k_n}\sum^{k_n}_{i=1}\log(X_{m_n,r+1}^{(i)})-a(m_n,r,p_n)\widehat{\gamma}_n\Big),
\]
where $\widehat{\gamma}_n$ is the estimator for $\gamma$ defined in \eqref{newest}, and $a(m,r,p)=\sum^m_{j=r+1}\frac1j+\log p$.

{   The estimator $\widehat{x}_p$ defined above is of the same
nature as the Weissman's  estimator for high quantiles based on a
full sample. In Weissman~\cite{Weissman1978}, the estimator for high
quantiles is a function of one extreme order statistic and an
estimator for the tail index, or more precisely, the estimator for
$\log x_p$ is the sum of the logarithm of one extreme order
statistic and the estimator of the tail index multiplied by a
constant. In our setting, we are able to define estimators for $\log
x_p$ by using the data from each of $k_n$ blocks, and thus we have
$k_n$ different estimators. Our estimator $\log \widehat{x}_p$ for
$\log x_p$ is the average of all these $k_n$ estimators. The only
difference is that we use a slightly different coefficient $a(m_n,
r, p)$ of $\widehat{\gamma}_n$ in our estimation in order to reduce
the bias in the estimation.  }

\begin{theorem}\label{thm2} In addition to conditions in Theorem~\ref{thm1},  if $p=p_n$ satisfies condition
 \[
m_np_n\to 0~\mbox{ and }~\log(m_np_n)=o(k_n^{1/2}) ~~\mbox { as }
n\to\infty,
 \]
then
\begin{equation}\label{qnorm}
\frac{(rk_n)^{1/2}}{|a(m_n,r,p_n)|}(\log\widehat{x}_p-\log x_p)\td
N(\delta b_r,\gamma^2),
\end{equation}
where $a(m,r,p)=\sum^m_{j=r+1}\frac1j+\log p$, and $b_r$ is defined
in \eqref{br}.
\end{theorem}

\vspace{10pt}

{   \noindent{\bf Remark 1.}  Since $k_n, m_n\to\infty$ and
$k_nm_m\sim n$ as $n\to\infty$, condition $m_np\to 0$ in
Theorem~\ref{thm2} is trivial when $np_n\to c\in [0,\infty)$. In
fact, we allow in Theorem~\ref{thm2} that $np_n\to \infty$ as long
as $m_np_n\to 0$ as $n\to\infty$. }

\vspace{10pt}

Now we consider the situation when the numbers of random variables
within the blocks are different and the numbers of the observations
available for inference are also different, and all random variables
are independent and identically distributed with a heavy-tailed
distribution (\ref{F}). Assume there are $k_n$ blocks of
observations, $V_i$, $1\le i\le k_n$, and in the $i$-th block $V_i$,
there are $m_i$ random variables, but only the $r_i+1$ largest order
statistics $X_{m_i, j}^{(i)}$, $j=1,\cdots, r_i+1$ are available for
inference. The total number of random variables within all $k_n$
blocks is $\sum^{k_n}_{i=1}m_i=n$ or $\sum^{k_n}_{i=1}m_i\sim n$.

Qi~\cite{Qi2010} proposed the following estimator for $\gamma$
\[
\widehat{\gamma}_n^*=
\frac{1}{\sum^{k_n}_{i=1}r_i}\sum^{k_n}_{i=1}\sum^{r_i}_{j=1}(\log
X_{m_i,j}^{(i)}-\log X_{m_i,r_i+1}^{(i)}).
\]

The asymptotic normality of $\widehat{\gamma}_n^*$ is obtained as
follows.

\begin{theorem}\label{thm3} (Qi~\cite{Qi2010})
 If (\ref{U})
holds and
 \[
k_n\to\infty, ~\frac{n}{k_n}\to \infty, \mbox{ and
}(\sum^{k_n}_{i=1}r_i)^{1/2}A(q_n)\to 0\mbox { as } n\to\infty,
 \]
 where $\displaystyle q_n=\min_{1\le i\le k_n}(m_i/r_i)\to \infty$ as $n\to\infty$, then
\[
(\sum^{k_n}_{i=1}r_i)^{1/2}(\widehat{\gamma}_n^*-\gamma){\stackrel{d}{\to}}
N(0,\gamma^2).
\]
\end{theorem}

\vspace{10pt}

Under the above setting-up,  we propose the following estimate for
$x_p$
\[
\widehat{x}^*_p=\exp\Big(\frac{1}{\sum^{k_n}_{j=1}r_j}\sum^{k_n}_{i=1}r_i\log(X_{m_i,r_i+1}^{(i)})-a_n(p_n)\widehat{\gamma}_n^*
\Big),
\]
where
\begin{equation}\label{anp}
a_n(p)=(\sum^{k_n}_{j=1}r_j)^{-1}\sum^{k_n}_{i=1}r_ia(m_i, r_i, p)
\end{equation}
with  $a(m, r, p)=\sum^m_{j=r+1}\frac1j+\log p$.

\begin{theorem}\label{thm4} In addition to conditions in Theorem~\ref{thm3}, assume the following condition
 \[
\max_{1\le i\le k_n}m_ip_n\to 0
 ~\mbox{ and }~ a_n(p_n)=o((\sum^{k_n}_{i=1}r_i)^{1/2}) ~~\mbox { as }
n\to\infty.
 \]
Then we have
\[
\frac{(\sum^{k_n}_{i=1}r_i)^{1/2}}{|a_n(p_n)|}(\log\widehat{x}^*_p-\log
x_p)\td N(0, \gamma^2).
\]
\end{theorem}

\vspace{10pt}

\noindent{\bf Remark 2.} Condition $\max_{1\le i\le k_n}m_ip_n\to 0$
is very weak. For example, when we consider a high quantile
$x_{p_n}$ under assumption that $np_n\to c\in [0,\infty)$, where
$n\sim\sum^{k_n}_{i=1}m_i$, condition $\max_{1\le i\le k_n}m_i/n\to
0$ implies $\max_{1\le i\le k_n}m_ip_n\to 0$. Under the conditions
in Theorem~\ref{thm4}, we can also show that $a_n(p_n)<0$
ultimately, see Lemma~\ref{lem2} in Appendix.


\vspace{10pt}

\noindent{\bf Remark 3.} In practice, handling the bias term in the
limiting distribution in Theorem~\ref{thm2} is not easy. One may
need to impose more restrictive conditions such as the so-called
third-order condition on $F$. We usually consider only the case
$\delta=0$ for convenience when we construct confidence intervals. A
$100(1-\alpha)\%$ confidence interval for $\log x_p$ based on the
normal approximation of $\log\widehat{x}_p$ in Theorem~\ref{thm2}
when $\delta=0$ is given by
\begin{equation}\label{CIN}
I_N(1-\alpha)=\left(\log\widehat{x}_p- z_{\alpha/2}\frac{|a(m,r,p)|{\widehat{\gamma}}_n}{(rk_n)^{1/2}},~
\log\widehat{x}_p+ z_{\alpha/2}\frac{|a(m,r,p)|{\widehat{\gamma}}_n}{(rk_n)^{1/2}}\right),
\end{equation}
where $z_{\alpha/2}$ is the critical value of the standard normal
distribution at level $\alpha/2$; that is,
$1-\Phi(z_{\alpha/2})=\alpha/2$, where $\Phi(x)$ is the cumulative
distribution function for the standard normal  random variable.
According to Theorem~\ref{thm2}, this confidence interval has an
asymptotically correct coverage probability, that is,
\[
P(\log x_{p,0}\in I_N(1-\alpha))\to 1-\alpha~~ \mbox{as}~~
n\to\infty,
\]
where $\log x_{p,0}$ is the true value of the parameter $\log x_p$.


\vspace{10pt}

 {   \noindent{\bf Remark 4.} The
second-order condition \eqref{U} is a standard condition that has
been used for univariate extreme-value statistics in the literature.
To verify condition \eqref{allkn}, the information on the function
$A$ is required. It is known that $A(t)$ is regularly varying at
infinity with index $\rho$. Theoretically, we can show that there
exists a function $g(t):=t^{-2\rho/(1-2\rho))}\ell(t)$, where
$\ell(t)$ is slowly varying at infinity, such that \eqref{allkn}
holds with $\delta=0$ if $k_n=o(g(n))$.  If a consistent estimator
of $\rho$, say $\hat\rho_n$, can be obtained, then for any fixed
$\varepsilon>0$, if
$k_n=o(n^{-2\hat\rho_n/(1-2\hat\rho_n)-\varepsilon})$, then
\eqref{allkn} holds with $\delta=0$.  When a complete sample is
available, consistent estimators for $\rho$ can be obtained, see,
e.g., Gomes et al.~\cite{GHP2002} and Peng and Qi~\cite{PQ2004}.
These estimators cannot be applied directly to the incomplete data
in the present paper.  Since the solution under our current setup
requires much effort, we leave this important work for the future
study. }

\vspace{10pt}


\section{Empirical likelihood and adjusted empirical likelihood methods}\label{elk}

In this section,  we assume $r\ge 1$ is a fixed integer, and $k_n$
and $m_n$ satisfy condition \eqref{allkn}. Set
\[
z_j^{(i)}(y)=j(\log X_{m,j}^{(i)}-\log X_{m,j+1}^{(i)})-\frac{1}{a(m,r,p)}(\log(X_{m,r+1}^{(i)})-y)
\]
for $j=1,\cdots, r$ and $i=1,\cdots,k_n$. Under conditions in
Theorem~\ref{thm2} with $\delta=0$, $\{z_j^{(i)}(y),~1\le j\le r,
1\le i\le k_n\}$ are approximately independent and identically
distributed with mean $0$ if $y=\log x_p$. We apply Owen's empirical
likelihood method (Owen~\cite{Owen1990}) to construct confidence
intervals or to test the hypothesis for logarithm of $x_p$.

 Let $\mathbf{q}=(q_1^{(1)}, \cdots,q_r^{(1)}, \cdots, q_1^{(k_n)},\cdots, q_r^{(k_n)})$
 be a probability vector
satisfying
\begin{equation}\label{probability}
\sum_{i=1}^{k_n}\sum^r_{j=1} q_j^{(i)}=1,  ~~q_j^{(i)} \geq 0~\mbox{ for } 1\le j\le r, ~1\le i\le k_n.
\end{equation}
Then the empirical likelihood, evaluated at  $y=\log x_p$, is
defined by
\begin{equation}\label{EL(y)}
EL(y)= \sup \left\{ \prod_{i=1}^{k_n}\prod^r_{j=1} q_j^{(i)}:
\sum_{i=1}^{k_n}\sum^r_{j=1} q_j^{(i)}z_j^{(i)}(y)=0, \sum_{i=1}^{k_n}\sum^r_{j=1} q_j^{(i)}=1\mbox{ with } q_j^{(i)}\ge 0\right\}.
\end{equation}
By the method of Lagrange multipliers, we can easily get the maximizers for the likelihood on the right-hand side of \eqref{EL(y)}
\[
q_j^{(i)}=\frac{1}{rk_n}\{ 1 + \lambda z_j^{(i)}(y)\}^{-1},
\;\; j=1,\cdots,r,~ i=1, \cdots, k_n,
\]
where $\lambda$ is the solution to the equation
\begin{equation}\label{eq40}
\sum_{i=1}^{k_n}\sum^r_{j=1}\frac{z_j^{(i)}(y)}{1 +
 \lambda z_j^{(i)}(y)}=0.
\end{equation}

On the other hand, $\prod_{i=1}^{k_n}\prod^r_{j=1} q_j^{(i)}$,
subject to constrains in \eqref{probability}, attains its maximum
$(rk_n)^{-rk_n}$ at $q_j^{(i)}=(rk_n)^{-1}$. So we define the
empirical likelihood ratio at $y_0$, the true value of $\log x_p$,
by
\[
 l(y_0)=\prod_{i=1}^{k_n}\prod^r_{j=1}(rk_nq_j^{(i)})=\prod_{i=1}^{k_n}\prod^r_{j=1}
  \{ 1 + \lambda z_j^{(i)}(y_0)\}^{-1},
\]
 and the corresponding empirical log-likelihood ratio statistic is defined as
\[
\mathcal{L}(y_0)=-2 \log l(y_0) = 2
\sum_{i=1}^{k_n}\sum^r_{j=1}\log
  \{ 1 + \lambda z_j^{(i)}(y_0)\},
\]
where $\lambda$ is the solution to (\ref{eq40}).

The following theorem gives the asymptotic distribution of
$\mathcal{L}(y_0)$.

\begin{theorem}\label{thm5} Under the conditions of Theorem~\ref{thm2} with $\delta=0$ we have
\begin{equation}\label{wilks}
 \mathcal{L}(y_0) {\stackrel{d}{\to}} \chi_1^2,
\end{equation}
where $\chi_1^2$ denotes a chi-squared random variable with one
degree of freedom, and $y_0$ is the true value of
$\log x_p$.
\end{theorem}

According to \eqref{wilks}, a $100(1-\alpha)\%$ confidence interval
for $\log x_p$ based on the empirical likelihood ratio statistic is
determined by
\[
I_E(1-\alpha)=\{y>0: \mathcal L(y)< c(\alpha)\},
\]
where $c(\alpha)$ is the $\alpha$ level critical value of a
chi-squared distribution with one degree of freedom.

When we define $EL(y)$ in \eqref{EL(y)}, we assume there is a
probability vector $\mathbf{q}$ satisfying \eqref{probability} such
that $\sum_{i=1}^{k_n}\sum^r_{j=1} q_j^{(i)}z_j^{(i)}(y)=0$, which
is equivalent to that $0$ is contained in the the convex hull of the
data points $\{z_j^{(i)}(y):~1\le j\le r, 1\le i\le k_n\}$. If this
is not true, the empirical likelihood ratio statistics $\mathcal
L(y_0)$ is set as infinity.  As a result, this may cause a serious
undercoverage for confidence intervals when the total number of data
points, $rk_n$, is relatively small. The same problem has been
discussed for the mean of iid random vectors in the literature, see,
e.g., Owen~\cite{Owen1988}, Hall and La Scala~\cite{HL1990}, Qin and
Lawless~\cite{QL1994}, and Tsao~\cite{Tsao2004}.

Chen et al.~\cite{CVA2008} proposed  the so-called  adjusted
empirical likelihood method to solve the undercoverage problem for
the mean of a distribution. For a random sample of size $n$, they
added a pseudo-sample point and applied the ordinary empirical
likelihood method to the $n+1$  data points. Later on,  Liu and
Chen~\cite{LC2010} investigated how to choose the correction factor
for the pseudo-sample point so as to achieve a better precision in
terms of coverage probability for empirical likelihood based
confidence intervals. Recent work by Li and Qi~\cite{LiQi2019}
applied the adjusted empirical likelihood method to constructing
confidence intervals for the tail index of a heavy-tailed
distribution.

Define a pseudo-data point as
\begin{equation}\label{zy}
z(y)=-\frac{a_n}{k_nr}\sum^{k_n}_{i=1}\sum^r_{j=1}z_j^{(i)}(y),
\end{equation}
where $a_n$ is a constant satisfying condition $a_n=o(k_n^{2/3})$.
In our applications, we will take $a_n=\frac{19}{12}$ which is the
optimal correction factor when the adjusted empirical likelihood is
applied to a random sample from an exponential distribution.  See
e.g., Li and Qi~\cite{LiQi2019} for more justifications.  The
so-called adjusted empirical likelihood method is to apply the
ordinary empirical likelihood method to the $k_nr+1$ data points
$\{z_j^{(i)}(y), 1\le i\le k_n, 1\le j\le r\} \cup\{z(y)\}$. By
following exactly the same procedure as the above, the adjusted
empirical likelihood ratio statistic at $y=\log x_{p}$ is given by
\[
\mathcal{AL}(y)= 2\sum_{i=1}^{k_n}\sum^r_{j=1}\log
  \{ 1 + \lambda z_j^{(i)}(y)\}+ 2\log
  \{ 1 + \lambda z(y)\},
\]
where $\lambda$ is the solution to the following equation
\[
 \sum_{i=1}^{k_n}\sum^r_{j=1}\frac{z_j^{(i)}(y)}{1 +
 \lambda z_j^{(i)}(y)}+\frac{z(y)}{1 +
 \lambda z(y)}  =0.
\]

\begin{theorem}\label{thm6} Assume the conditions of Theorem~\ref{thm2} with $\delta=0$ are satisfied
and $a_n=o(k_n^{2/3})$ as $n\to\infty$. Then we have
\[
\mathcal{AL}(y_0) {\stackrel{d}{\to}} \chi_1^2,
\]
where $\chi_1^2$ denotes a chi-squared random variable with one
degree of freedom, and $y_0$ is the true value of
$\log x_p$.
\end{theorem}

According to Theorem~\ref{thm6}, a $100(1-\alpha)\%$ confidence interval
for $\log x_p$ based on the adjusted  empirical likelihood ratio statistic is
determined by
\begin{equation}\label{AECIN}
I_{AE}(1-\alpha)=\{y>0: \mathcal{AL}(y)< c(\alpha)\},
\end{equation}
where $c(\alpha)$ is the $\alpha$ level critical value of a
chi-squared distribution with one degree of freedom.

\section{Simulation study}\label{simulation}

In this section, we carry out a simulation study to compare the
performance of the confidence intervals based on the adjusted
empirical likelihood ($I_{AE}^*(1-\alpha)$ defined in \eqref{AECIN})
and the normal approximation ($I_N(1-\alpha)$ given in \eqref{CIN})
for high quantiles in terms of coverage probability and interval
length. We take $a_n=\frac{19}{12}$ for the weight $a_n$ in
\eqref{zy}. We consider the following two types of cumulative
distribution functions (cdf):
\begin{description}
  \item[(a).] the Fr\'{e}chet cdf given by $F(x)=\exp(-x^{-a})$ $(x>0),$ where
$a>0$ (notation: Fr\'echet($a$));
  \item[(b).] the Burr cdf
given by $F(x)=1-(1+x^{a})^{-b}$ $(x>0),$ where $a>0$, $b>0$
(notation: Burr($a,b$)).
\end{description}

In our simulation study, we choose $r=1$, that is, we consider the
case when only two largest observations within blocks are used for
the inference.
 We choose the confidence level $1-\alpha=95\%$ in the study.   The
simulation is implemented by Software $\mathbf{R}$. We will use the
following three distributions in our study:  Fr\'echet ($1$), Burr
($0.5, 1$) and Burr ($1,0.5$). Both Fr\'echet and Burr distributions
can be expanded in the form given in \eqref{2terms}, and their
second-order function $A(t)$  is proportional to $t^{\rho}$, where
$\rho=-1$ for Fr\'echet ($1$) and Burr ($0.5, 1$) and $\rho=-2$ for
Burr($1,0.5$).


From each of three distributions, Fr\'echet($1$), Burr($0.5, 1$) and
Burr($1, 0.5$), we generate $k$ blocks of independent random
variables with $k=10, 15, \cdots, 100$, and each block contains $m$
observations, where $m$ can be selected under one of the following
two schemes. For each distribution and each combination of $k$
 and $m$,  the coverage probabilities and average lengths of two confidence
 intervals, $I_E^*(0.95)$ and $I_N(0.95)$), are estimated based on $5000$ replicates.


 \noindent\underline{\it Scheme 1}. Set $n=1000$ and $p=p_n=1/n$.
 For each $k$ from $\{10, 15, \cdots, 100\}$, define $m=[1000/k]$,
where $[x]$ denotes the integer part of $x$.  This setup applies to
all three distributions. The coverage probabilities and the average
lengths of confidence intervals based on the adjusted empirical
likelihood method (AELM) and normal approximation method (NORM) for
$\log x_p$ are obtained and reported in Tables~\ref{table1} and
\ref{table2}, respectively.

\noindent\underline{\it Scheme 2}. For each $k$ from $\{10, 15,
\cdots, 100\}$, we select $m$ as a function $k$ in the form
$m=[50k^{v}]$, where $v\in (0,1)$ depends on the distribution from
which the observations are sampled. In particular, we have selected
$v=1/2$ for Burr($0.5,1$), $v=1/4$ for Burr($1,0.5$), and $v=1/2$
for Fr\'{e}chet(1). Notice that $m$ increases with $k$, and so does
the total number, $km$, of observations within all $k$ blocks. For
each combination of $k$ and $m$, we estimate $\log x_p$ with
$p=1/(km)$. Again, we estimate the coverage probability and average
length of two confidence intervals, $I_E^*(0.95)$ based on the
adjusted empirical likelihood method (AELM) and $I_N(0.95)$ based
normal approximation method (NORM). Simulation results are given in
Tables~\ref{table3} and \ref{table4}.

Under Scheme 2,  with the specific selection $m=[50k^{v}]$ for each
distribution, $k^{1/2}A(m)$ is approximatively a constant as $k$
gets larger, and the multiplier $50$ is selected so that the bias
term of the limiting distribution in \eqref{qnorm} is relatively
small. Under Scheme 1,  the total number of observations in the $k$
blocks is approximately $1000$, and $m$ decreases with $k$, for
example,  $m=100$ if $k=10$, and $m=10$ if $k=100$. Since $|A(t)|$
is proportional to $t^{\rho}$ for some $\rho<0$, we see that the
bias term in \eqref{qnorm} is getting larger when $k$ is bigger
under Scheme 1.

From Tables~\ref{table1} and \ref{table3}, the confidence intervals
based on the normal approximation have a significantly lower
coverage for smaller values of $k$, and those based on the adjusted
empirical likelihood method achieve a much better coverage in this
case. The performance of the normal approximation is getting better
when $k$ increases under Scheme 1. For Burr($1,0.5$), the
performance of the adjusted empirical likelihood method may be
greatly influenced by the bias terms which are approximately
proportional to $k^{1/2}m^{-2}\sim k^{2.5}/1000^2$. Under Scheme 2,
$m$ increases with $k$ and the bias terms in \eqref{qnorm} are
reasonably small, and the coverage probabilities from
Table~\ref{table3} indicate that the adjusted empirical likelihood
method outperforms the normal approximation method for all three
distributions.

Tables~\ref{table2} and \ref{table4} reveal that the average lengths
of the confidence intervals based on the normal approximation are
shorter than those based on the adjusted empirical likelihood method
under both Scheme 1 and Scheme 2 when $k$ is relatively small. This
is caused by the lower coverage of confidence intervals based on the
normal approximation.  When $k$ is getting larger, the average
lengths for both methods are comparable.

In summary, we conclude that the adjusted empirical likelihood
method results in better coverage probability when $k$ is small or
when the bias term in \eqref{qnorm} is relatively small. When the
bias term in \eqref{qnorm} is getting too large, theoretically, both
methods have an undercoverage problem, but the adjusted empirical
likelihood method may suffer more than the normal approximation
method as we have observed in Table~\ref{table1} for Burr($1,0.5$)
distribution.

{   Finally, we conclude this section with some comments on
application of the proposed methods in the paper.  In general, a
relatively large sample is required for application of results under
the framework of extreme value statistics since only a few of
largest observations in the sample can be used in the estimation.
The accuracy of the normal approximation for estimators of $\gamma$
and $x_p$ depends on the number of observations used in the
estimation ($rk_n$) and the asymptotic bias for the normalized
statistics $rk_n^{1/2}A(n/k_n)$. We recommend that $rk_n\ge 30$ and
$rk_n^{1/2}|A(n/k_n)|\le 0.2$.  Since $A(t)$ is proportional to
$t^{\rho}$ in most applications,  we can assume
$rk_n^{1/2}(n/k_n)^{\rho}\le 0.2$.  Now consider the special case
$r=1$, we have $n\ge 5^{-1/\rho}k_n^{(1-2\rho)/(-2\rho)}$ as
suggested sample size for $k_n\ge 30$.  Since $\rho$ is unknown,
estimation of $\rho$ seems important issue for this purpose. See
more comments in Remark 4 at the end of Section~\ref{main}.

}

\begin{table}[!htb]
\caption{Coverage probabilities for adjusted empirical likelihood
method with correction factor $19/12$ (AELM) and normal
approximation method for $\log\widehat{x}_p$ (NORM): the number of
observations within each block is set to be  $m=[1000/k$] (Scheme
1)}
\label{table1}\centering {
 \footnotesize
\begin{tabular}{r|rr|rr|rr}\hline
  \hline
    & \multicolumn{2}{c|}{Fr\'echet($1$)} &\multicolumn{2}{c|}{Burr($0.5$,$1$)}&\multicolumn{2}{c}{Burr($1$,$0.5$)}\\
   $k_n$ & AELM & NORM& AELM & NORM& AELM & NORM\\
\hline
 10&  0.9630 &   0.8996 & 0.9602 &   0.9046 & 0.9612  &  0.9066\\
 15&  0.9420 &   0.9172 & 0.9342 &   0.9114 & 0.9354  &  0.9186\\
 20&  0.9372 &   0.9256 & 0.9360 &   0.9234 & 0.9384  &  0.9286\\
 25&  0.9408 &   0.9294 & 0.9410 &   0.9308 & 0.9438  &  0.9394\\
 30&  0.9440 &   0.9364 & 0.9406 &   0.9248 & 0.9448  &  0.9348\\
 35&  0.9438 &   0.9412 & 0.9494 &   0.9388 & 0.9524  &  0.9520\\
 40&  0.9448 &   0.9490 & 0.9442 &   0.9384 & 0.9434  &  0.9522\\
 45&  0.9490 &   0.9498 & 0.9430 &   0.9370 & 0.9440  &  0.9506\\
 50&  0.9490 &   0.9510 & 0.9462 &   0.9446 & 0.9440  &  0.9582\\
 55&  0.9446 &   0.9482 & 0.9374 &   0.9358 & 0.9364  &  0.9488\\
 60&  0.9484 &   0.9534 & 0.9460 &   0.9418 & 0.9418  &  0.9574\\
 65&  0.9498 &   0.9600 & 0.9470 &   0.9452 & 0.9414  &  0.9576\\
 70&  0.9464 &   0.9566 & 0.9488 &   0.9438 & 0.9348  &  0.9592\\
 75&  0.9494 &   0.9610 & 0.9470 &   0.9472 & 0.9300  &  0.9602\\
 80&  0.9458 &   0.9572 & 0.9420 &   0.9434 & 0.9258  &  0.9548\\
 85&  0.9436 &   0.9570 & 0.9458 &   0.9454 & 0.9146  &  0.9534\\
 90&  0.9498 &   0.9616 & 0.9446 &   0.9464 & 0.9192  &  0.9510\\
 95&  0.9408 &   0.9580 & 0.9468 &   0.9490 & 0.9044  &  0.9492\\
100&  0.9384 &   0.9538 & 0.9462 &   0.9502 & 0.8988  &  0.9438\\
\hline
\end{tabular}}
\end{table}

\begin{table}[!htb]
\caption{Averages of lengths for confidence intervals based adjusted
empirical
 likelihood method with correction factor $19/12$ (AELM) and normal approximation method for $\log\widehat{x}_p$
 (NORM): the number of
observations within each block is set to be  $m=[1000/k$] (Scheme
1)}
\label{table2}\centering {
\footnotesize
\begin{tabular}{r|rr|rr|rr}\hline
  \hline
    & \multicolumn{2}{c|}{Fr\'echet($1$)} &\multicolumn{2}{c|}{Burr($0.5$,$1$)}&\multicolumn{2}{c}{Burr($1$,$0.5$)}\\
   $k_n$ & AELM & NORM& AELM & NORM & AELM & NORM\\
\hline
 10& 5.014 & 3.393 &9.985 & 6.754 &10.052 & 6.834 \\
 15& 3.622 & 3.193 &7.163 & 6.343 & 7.184 & 6.373 \\
 20& 3.284 & 3.015 &6.525 & 5.979 & 6.550 & 6.027 \\
 25& 3.082 & 2.875 &6.159 & 5.702 & 6.170 & 5.761 \\
 30& 2.946 & 2.780 &5.846 & 5.470 & 5.854 & 5.526 \\
 35& 2.818 & 2.683 &5.605 & 5.276 & 5.608 & 5.317 \\
 40& 2.702 & 2.588 &5.386 & 5.095 & 5.392 & 5.126 \\
 45& 2.623 & 2.524 &5.199 & 4.937 & 5.216 & 4.979 \\
 50& 2.537 & 2.451 &5.040 & 4.801 & 5.060 & 4.841 \\
 55& 2.471 & 2.401 &4.895 & 4.682 & 4.909 & 4.707 \\
 60& 2.429 & 2.369 &4.806 & 4.609 & 4.784 & 4.594 \\
 65& 2.363 & 2.315 &4.670 & 4.494 & 4.672 & 4.489 \\
 70& 2.307 & 2.265 &4.560 & 4.396 & 4.564 & 4.399 \\
 75& 2.263 & 2.228 &4.468 & 4.321 & 4.466 & 4.316 \\
 80& 2.228 & 2.201 &4.393 & 4.259 & 4.380 & 4.233 \\
 85& 2.199 & 2.181 &4.322 & 4.204 & 4.294 & 4.160 \\
 90& 2.135 & 2.119 &4.198 & 4.086 & 4.219 & 4.092 \\
 95& 2.118 & 2.111 &4.165 & 4.070 & 4.145 & 4.023 \\
100& 2.062 & 2.057 &4.056 & 3.967 & 4.079 & 3.962 \\
 \hline
\end{tabular}}
\end{table}


\begin{table}[!htb]
\caption{Coverage probabilities for adjusted empirical likelihood
method with correction factor $19/12$ (AELM) and normal
approximation method for $\log\widehat{x}_p$ (NORM): the number of
observations within each block is set to be  $m=[50k^v]$ (Scheme 2),
where $v=1/2$ for Burr($0.5,1$), $v=1/4$ for Burr($1,0.5$), and
$v=1/2$ for Fr\'{e}chet(1)}
\label{table3}\centering {
\footnotesize
\begin{tabular}{r|rr|rr|rr}\hline
  \hline
    & \multicolumn{2}{c|}{Fr\'echet($1$)} &\multicolumn{2}{c|}{Burr($0.5$,$1$)}&\multicolumn{2}{c}{Burr($1$,$0.5$)}\\
   $k_n$ & AELM & NORM& AELM & NORM& AELM & NORM\\
\hline
 10&   0.9648 &   0.9020& 0.9578&    0.8966& 0.9604 &   0.9036 \\
 15&   0.9416 &   0.9180& 0.9378&    0.9148& 0.9402 &   0.9166 \\
 20&   0.9366 &   0.9216& 0.9388&    0.9248& 0.9370 &   0.9216 \\
 25&   0.9422 &   0.9286& 0.9360&    0.9294& 0.9398 &   0.9312 \\
 30&   0.9422 &   0.9284& 0.9386&    0.9306& 0.9410 &   0.9356 \\
 35&   0.9410 &   0.9350& 0.9388&    0.9294& 0.9392 &   0.9340 \\
 40&   0.9456 &   0.9366& 0.9418&    0.9332& 0.9406 &   0.9384 \\
 45&   0.9470 &   0.9392& 0.9434&    0.9368& 0.9462 &   0.9440 \\
 50&   0.9472 &   0.9390& 0.9436&    0.9346& 0.9462 &   0.9386 \\
 55&   0.9440 &   0.9356& 0.9460&    0.9394& 0.9422 &   0.9362 \\
 60&   0.9464 &   0.9398& 0.9456&    0.9324& 0.9448 &   0.9414 \\
 65&   0.9480 &   0.9418& 0.9434&    0.9322& 0.9416 &   0.9356 \\
 70&   0.9492 &   0.9448& 0.9458&    0.9390& 0.9456 &   0.9402 \\
 75&   0.9464 &   0.9418& 0.9460&    0.9406& 0.9480 &   0.9440 \\
 80&   0.9474 &   0.9452& 0.9474&    0.9432& 0.9510 &   0.9480 \\
 85&   0.9536 &   0.9448& 0.9488&    0.9396& 0.9486 &   0.9448 \\
 90&   0.9448 &   0.9434& 0.9442&    0.9356& 0.9524 &   0.9508 \\
 95&   0.9500 &   0.9436& 0.9508&    0.9458& 0.9518 &   0.9484 \\
100&   0.9464 &   0.9374& 0.9468&    0.9402& 0.9478 &   0.9440 \\
\hline
\end{tabular}}
\end{table}

\begin{table}[!htb]
\caption{Averages of lengths for confidence intervals based adjusted
empirical
 likelihood method with correction factor $19/12$ (AELM) and normal approximation method for $\log\widehat{x}_p$
 (NORM):  the number of
observations within each block is set to be  $m=[50k^v]$ (Scheme 2),
where $v=1/2$ for Burr($0.5,1$), $v=1/4$ for Burr($1,0.5$), and
$v=1/2$ for Fr\'{e}chet(1)}
\label{table4}\centering {
\footnotesize
\begin{tabular}{r|rr|rr|rr}\hline
  \hline
    & \multicolumn{2}{c|}{Fr\'echet($1$)} &\multicolumn{2}{c|}{Burr($0.5$,$1$)}&\multicolumn{2}{c}{Burr($1$,$0.5$)}\\
   $k_n$ & AELM & NORM& AELM & NORM & AELM & NORM\\
\hline
 10& 5.008 & 3.386 &9.933 & 6.701& 10.052 & 6.834 \\
 15& 3.591 & 3.172 &7.160 & 6.311&  7.184 & 6.373 \\
 20& 3.280 & 2.994 &6.529 & 5.986&  6.550 & 6.027 \\
 25& 3.078 & 2.856 &6.141 & 5.700&  6.170 & 5.761 \\
 30& 2.919 & 2.733 &5.832 & 5.469&  5.854 & 5.526 \\
 35& 2.798 & 2.637 &5.586 & 5.271&  5.608 & 5.317 \\
 40& 2.690 & 2.547 &5.376 & 5.090&  5.392 & 5.126 \\
 45& 2.603 & 2.475 &5.186 & 4.935&  5.216 & 4.979 \\
 50& 2.524 & 2.404 &5.032 & 4.800&  5.060 & 4.841 \\
 55& 2.452 & 2.343 &4.895 & 4.676&  4.909 & 4.707 \\
 60& 2.395 & 2.293 &4.768 & 4.557&  4.784 & 4.594 \\
 65& 2.333 & 2.238 &4.652 & 4.463&  4.672 & 4.489 \\
 70& 2.280 & 2.192 &4.544 & 4.375&  4.564 & 4.399 \\
 75& 2.235 & 2.149 &4.449 & 4.287&  4.466 & 4.316 \\
 80& 2.190 & 2.109 &4.356 & 4.204&  4.380 & 4.233 \\
 85& 2.150 & 2.073 &4.274 & 4.127&  4.294 & 4.160 \\
 90& 2.110 & 2.036 &4.201 & 4.057&  4.219 & 4.092 \\
 95& 2.071 & 2.001 &4.126 & 3.989&  4.145 & 4.023 \\
100& 2.037 & 1.970 &4.059 & 3.936&  4.079 & 3.962 \\
 \hline
\end{tabular}}
\end{table}

\section{Proofs}\label{proofs}

Before we prove the main results, we introduce some notations. As we
have assumed that $X_i\ge 1$ for $i\ge 1$,  we see that $F^-(u)\ge
1$ for all $u\in (0,1)$, and thus $U(x)=F^-(1-\frac1x)\ge 1$ is well
defined for $x>1$. Note that  $U(x)$ is non-decreasing for $x>1$.

From de Haan and Stadtm\"{u}ller~\cite{deHannS1996},  condition
(\ref{U}) implies that $A(t)$ is a regularly varying function with
index $\rho$, and $|A(t)|$ is bounded away from $0$ and $\infty$ on
every compact subset of $[c_0,\infty)$, where $c_0>0$ is a constant.
Without loss of generality, we assume $A(t)$ is bounded away from 0
and $\infty$ in $(0, c_0]$ since redefining $A(t)$ on $(0, c_0)$
doesn't change condition (\ref{U}).  Then using Potter's bounds to
$A(x)$, for every $\delta>0$, there exists a constant $c(\delta)>0$
such that
\begin{equation}\label{abound}
|\frac{A(x)}{A(y)}|\le
c(\delta)\max\big((\frac{x}{y})^{\rho+\delta},
(\frac{x}{y})^{\rho-\delta}\big)~~\mbox{ for all }x,y>0.
\end{equation}
See, e.g., Theorem 1.5.6 in  Bingham {\it et al.}~\cite{BGT1987}.

Next, we see that (\ref{U}) is equivalent to
\begin{equation}\label{logu}
\lim_{t\to\infty}\frac{\log U(tx)-\log U(t)-\gamma \log
x}{A(t)}=\frac{x^{\rho}-1}{\rho},~~x>0.
 \end{equation}


Let $h(x)=\log U(x)-\gamma \log x$, $x>1$. Then (\ref{logu}) implies that for each
$x>0$,
\[
\frac{h(tx)-h(t)}{A(t)}=\frac{\log
U(tx)-\log U(t)-\gamma \log x}{A(t)}\to \frac{x^{\rho}-1}{\rho}~\mbox{ as }t\to\infty.
\]
Since $\rho<0$, we have from Theorem B.2.18 in de Haan and
Ferreira~\cite{deHannF2006} that there exists a constant $c$ such
that $\lim_{x\to\infty}h(x)=c$, and $h_1(x):=h(t)-c$ is a regularly
varying function with index $\rho$, i.e.,
\[
\frac{h_1(tx)-h_1(t)}{A(t)}=\frac{\log
U(tx)-\log U(t)-\gamma \log x}{A(t)}\to \frac{x^{\rho}-1}{\rho}~\mbox{ as }t\to\infty.
\]
In this case, we have  $A(x)\sim \rho h_1(x)$ as $x\to\infty$,
$h_1(x)\to 0$ as $x\to\infty$, and $h_1(x)$ is uniformly bounded in
interval $[1, \infty)$.  Meanwhile, we have $|h_1(y)/A(y)|$ is
uniformly bounded in $(1, \infty)$, and we assume $|h_1(y)/A(y)|\le
C_0$ for some $C_0>0$.

Rewrite
\[
\frac{\log U(tx)-\log U(t)-\gamma \log
x}{A(t)}=\frac{h_1(tx)}{A(tx)}\frac{A(tx)}{A(t)}-\frac{h_1(t)}{A(t)},~~t>0,
tx>1.
\]
By substituting $y$ for $x$ in the above equation and subtracting it from the above equation we have
\begin{eqnarray*}
|\frac{\log U(tx)-\log U(ty)-\gamma(\log x-\log y)}{A(t)}|&=&|\frac{h_1(tx)}{A(tx)}\frac{A(tx)}{A(t)}-\frac{h_1(ty)}{A(ty)}\frac{A(ty)}{A(t)}|\\
&\le& C_0(|\frac{A(tx)}{A(t)}|+|\frac{A(ty)}{A(t)}|).
\end{eqnarray*}
Now we apply Potter's bounds \eqref{abound} to both
$\frac{A(tx)}{A(t)}$ and $\frac{A(ty)}{A(t)}$ with $\delta=-\rho/2$
and conclude that
\begin{equation}\label{logbound}
|\frac{\log U(tx)-\log U(ty)-\gamma(\log x-\log y)}{A(t)}|\le
C_1(x^{\rho/2}+x^{3\rho/2}+y^{\rho/2}+y^{3\rho/2})
\end{equation}
for all $t>1$, $tx>1$, and $ty>1$,  where $C_1>0$ is a constant.


As in Qi~\cite{Qi2010}, our proofs rely on the distributional
representations for the observations.  We will use the same notation
as in Qi~\cite{Qi2010}.

Assume $\{E_j^{(i)},~i, j\ge 1\}$ are iid random variables with a
unit exponential distribution. It is easy to see that
$\{U(e^{E_j^{(i)}}), ~i, j\ge 1\}$ are iid random variables with the
distribution $F$.

Apparently, $\{X_{m_i,j}^{(i)}, 1\le j\le m_i\}$ have the same joint
distribution as $\{U(E_{m_i,j}^{(i)}), 1\le j\le m_i\}$, where
$E_{m_i,1}^{(i)}\ge \cdots\ge E_{m_i, m_i}^{(i)}$ are the order
statistics of $E_{j}^{(i)}, 1\le j\le m_i$. Without loss of
generality, we assume that
\begin{equation}\label{RepX}
X_{m_i,j}^{(i)}=U(e^{E_{m_i,j}^{(i)}}), ~~~ 1\le j\le m_i,~1\le i\le
k_n.
\end{equation}

For each $i\ge 1$, set
$I_j^{(i)}=j(E_{m_i,j}^{(i)}-E_{m_i,j+1}^{(i)})$ for $j=1,\cdots,
m_i-1$ and $I_{m_i}^{(i)}=m_iE_{m_i,m_i}^{(i)}$. Then
$\{I_{j}^{(i)}, ~1\le j\le m_i, ~1\le i\le k_n\}$ are iid random
variables with a unit exponential distribution. We also have
\begin{equation}\label{Exrep}
E_{m_i, r+1}^{(i)}=\sum^{m_i}_{j=r+1}\frac{I_j^{(i)}}{j}.
\end{equation}
It is easy to see that $E_{m_i,r+1}^{(i)}$, $i\ge 1$ are independent
random variables with their means and variances given by
\begin{equation}\label{mean-var}
\E(E_{m_i,r+1}^{(i)})=\sum^{m_i}_{j=r+1}\frac1j,~~\mathrm{Var}(E_{m_i,r+1}^{(i)})=\sum^{m_i}_{j=r+1}\frac1{j^2}\le\frac1r.
\end{equation}


\begin{lemma}\label{lem-sum} As $n\to\infty$,
\[
\frac1{(\sum^{k_n}_{i=1}r_i)^{1/2}}\sum^{k_n}_{i=1}r_i(E_{m_i,r_i+1}^{(i)}-\sum^{m_n}_{j=r_i+1}\frac1j)=O_p(1).
\]
\end{lemma}

\noindent{\it Proof.} The lemma is trivial since the variance or the
second moment of the left-hand side above is equal to
\[
\frac1{\sum^{k_n}_{i=1}r_i}\sum^{k_n}_{i=1}r_i^2\sum^{m_i}_{j=r_i+1}\frac1{j^2}\le
\frac1{\sum^{k_n}_{i=1}r_i}\sum^{k_n}_{i=1}r_i^2\frac{1}{r_i}=1
\]
from \eqref{Exrep} and \eqref{mean-var}.
 \qed

\begin{lemma}\label{lem1}  For any $\delta>0$ we have
\[
\E\Big((\frac{\exp(E_{m, r+1}^{(1)})}{m/r})^{-\delta}\Big)\le
\exp(\delta+\frac{\delta^2}{2}),~~1\le r<m,~m\ge 2.
\]
\end{lemma}

\noindent{\it Proof.}  Using representation \eqref{Exrep} and the moment-generating function of exponential random variables we have
\begin{eqnarray*}
\E\big((\frac{\exp(E_{m,r+1}^{(1)})}{m/r})^{-\delta}\big)&=&(m/r)^{\delta}\E(\exp(-\delta \sum^m_{j=r+1}\frac{I_j^{(1)}}{j} ))\\
&=&\frac{(m/r)^{\delta}}{\prod^m_{j=r+1}(1+\frac{\delta}{j})}\\
&=&\exp\big(\delta\log(m/r)-\sum^m_{j=r+1}\log(1+\frac{\delta}{j})\big)\\
&\le &\exp\big(\delta\log (m/r)-\sum^m_{j=r+1}(\frac{\delta}{j}-\frac{\delta^2}{2j^2})\big)\\
&\le &\exp\big(\delta(\log (m/r)-\sum^m_{j=r+1}\frac{1}{j})+\frac{\delta^2}{2}\sum^m_{j=r+1}{j^{-2}})\big)\\
&\le&\exp(\delta+\frac{\delta^2}{2}).
\end{eqnarray*}
In the above estimation we have used inequalities that $\log(1+y)\ge
y-\frac12y^2$ for $y>0$ and $\log
(m/r)-\sum^m_{j=r+1}\frac{1}{j}<\frac{1}{r}\le 1$. \qed

\begin{lemma}\label{lem2} Under conditions $\min_{1\le i\le k_n}(m_i/r_i)\to \infty$ and $\max_{1\le i\le k_n}m_ip_n\to
0$, we have $\min_{1\le i\le k_n}(-a(m_i, r_i, p_n))\to\infty$ and
$-a_n(p_n)\to\infty$ as $n\to\infty$.
\end{lemma}

\noindent{\it Proof.}  Since
\[
\sum^m_{j=r+1}\frac{1}{j}<\int^m_{r}\frac1xdx=\log(\frac{m}{r})<\sum^m_{j=r}\frac{1}{j}=\frac{1}{r}+\sum^m_{j=r+1}\frac{1}{j}
\]
for $1\le r<m$, we have
\[
\log(\frac{m}{r})-\frac{1}{r}<\sum^m_{j=r+1}\frac{1}{j}<\log(\frac{m}{r}),
\]
which implies
\[
\log(\frac{mp_n}{r})-\frac{1}{r}<\sum^m_{j=r+1}\frac{1}{j}+\log
p_n<\log(\frac{mp_n}{r}).
\]
Therefore, for $1\le i\le k_n$,
\[
\log(\frac{r_i}{m_ip_n})<-a(m_i,r_i,
p_n)<\log(\frac{r_i}{m_ip_n})+\frac{1}{r_i}.
\]
Since $\min_{1\le i\le k_n}\frac{r_i}{m_ip_n}\ge \frac{1}{\max_{1\le
i\le k_n}m_ip_n}\to\infty$, we obtain
\[
\min_{1\le i\le k_n}(-a(m_i,r_i, p_n))\sim \min_{1\le i\le
k_n}\log(\frac{r_i}{m_ip_n})\to\infty
\]
as $n\to\infty$. This also implies $-a_n(p_n)\to\infty$ from
definition \eqref{anp}. \qed

We will prove a general result which can be used in the proofs for
Theorems~\ref{thm2} and \ref{thm4}.

\begin{lemma}\label{lem3} Under conditions $q_n=\min_{1\le i\le k_n}(m_i/r_i)\to\infty$ and $\max_{1\le i\le k_n}m_ip_n\to 0$,  we have
\[
\frac{1}{\sum^{k_n}_{i=1}r_i}\sum^{k_n}_{i=1}r_i\log(X_{m_i,r_i+1}^{(i)})-\log
x_p=\frac{\sum^{k_n}_{i=1}r_ia(m_i,r_i,p)}{\sum^{k_n}_{i=1}r_i}\gamma
+O_p(\frac{1}{(\sum^{k_n}_{i=1}r_i)^{1/2}}+|A(q_n)|).
\]
\end{lemma}

\noindent{\it Proof.} Write
\[
\varepsilon(t,x,y)=\frac{\log U(tx)-\log U(ty)-\gamma(\log x-\log y)}{A(t)}.
\]
Then  from \eqref{logbound} we have
\begin{equation}\label{error}
|\varepsilon(t,x,y)|\le  C_1(x^{\rho/2}+x^{3\rho/2}+y^{\rho/2}+y^{3\rho/2})
\end{equation}
for $t>1$, $tx>1$, $ty>1$, and
\[
\log U(tx)-\log U(ty)=\gamma(\log x-\log y)+A(t)\varepsilon(t,x,y).
\]

Review that $x_p=U(\frac{1}{p_n})$.  For each $i\ge 1$, by using
representation \eqref{RepX} with $t=m_i/r_i$,
$tx=e^{E_{m_i,r_i+1}^{(i)}}$ and $ty=\frac{1}{p_n}$ we have
\begin{eqnarray*}
\log(X_{m_i,r_i+1}^{(i)})-\log x_p&=&\log U(e^{E_{m_i,r_i+1}^{(i)}})-\log U(\frac1{p})\\
&=&\gamma(\log\frac{e^{E_{m_i,r_i+1}^{(i)}}}{m_i/r_i}-\log\frac{r_i}{m_ip_n})+A(\frac{m_i}{r_i})\varepsilon(
\frac{m_i}{r_i},\frac{e^{E_{m_i,r_i+1}^{(i)}}}{m_i/r_i},\frac{r_i}{m_ip_n})\\
&=&\gamma(E_{m_i,r_i+1}^{(i)}+\log p_n)
+A(\frac{m_i}{r_i})\varepsilon(\frac{m_i}{r_i},\frac{e^{E_{m_i,r_i+1}^{(i)}}}{m_i/r_i},\frac{r_i}{m_ip_n}).
\end{eqnarray*}
Then
\begin{eqnarray}\label{RepXPP}
&&\log(X_{m_i,r_i+1}^{(i)})-\log x_p- \gamma
a(m_i,r_i,p_n)\nonumber\\
&=&\gamma(E_{m_i,r_i+1}^{(i)}-\sum^m_{j=r_i+1}\frac1j)
+A(\frac{m_i}{r_i})\varepsilon(\frac{m_i}{r_i},\frac{e^{E_{m_i,r_i+1}^{(i)}}}{m_i/r_i},\frac{r_i}{m_ip_n}).
\end{eqnarray}
We have from \eqref{error} that
\[
|\varepsilon(\frac{m_i}{r_i},\frac{e^{E_{m_i,r_i+1}^{(i)}}}{m_i/r_i},\frac{r_i}{m_ip_n})|\le
C_1\Big(\big(\frac{e^{E_{m_i,r_i+1}^{(i)}}}{m_i/r_i}\big)^{\rho/2}+
\big(\frac{e^{E_{m_i,r_i+1}^{(i)}}}{m_i/r_i}\big)^{3\rho/2}+2\Big)
\]
as long as $\frac{m_i}{r_i}>1$ and $\frac{r_i}{m_ip_n}>1$, which are
true for all $1\le i\le k_n$ since  $\max_{1\le i\le k_n}m_ip_n\to
0$ and $q_n=\min_{1\le i\le k_n}(m_i/r_i)\to\infty$ as $n\to\infty$.

From Lemma~\ref{lem1}, for any $d>0$,
$\E\big(|\varepsilon(\frac{m_i}{r_i},\frac{e^{E_{m_i,r_i+1}^{(i)}}}{m_i/r_i},\frac{r_i}{m_ip_n})|^d\big)$
are uniformly bounded for $1\le i\le k_n$ for all large $n$. We can
conclude that
\begin{equation}\label{maxep}
\max_{1\le i\le
k_n}|\varepsilon(\frac{m_i}{r_i},\frac{e^{E_{m_i,r_i+1}^{(i)}}}{m_i/r_i},\frac{r_i}{m_ip_n})|=O_p(k_n^{1/2})
\end{equation}
and
\begin{equation}\label{sumep}
\frac{1}{\sum^{k_n}_{i=1}r_i}\sum^{k_n}_{i=1}r_i|\varepsilon(\frac{m_i}{r_i},\frac{e^{E_{m_i,r_i+1}^{(i)}}}{m_i/r_i},\frac{r_i}{m_ip_n})|^\ell=O_p(1)
\end{equation}
for any $\ell=1,2$.
\eqref{maxep} is true since there exists a $C>0$ such that for any $x>0$
\begin{eqnarray*}
&&P(\max_{1\le i\le k_n}|\varepsilon(\frac{m_i}{r_i},\frac{e^{E_{m_i,r_i+1}^{(i)}}}{m_i/r_i},\frac{r_i}{m_ip_n})|>xk_n^{1/2})\\
&\le&  \sum^{k_n}_{i=1}P(|\varepsilon(\frac{m_i}{r_i},\frac{e^{E_{m_i,r_i+1}^{(i)}}}{m_i/r_i},\frac{r_i}{m_ip_n})|>xk_n^{1/2})\\
&\le& \sum^{k_n}_{i=1}\frac{\E(\varepsilon(\frac{m_i}{r_i},\frac{e^{E_{m_i,r_i+1}^{(i)}}}{m_i/r_i},\frac{r_i}{m_ip_n}))^2}{x^2k_n}\\
&\le &\frac{C}{x^2}.
\end{eqnarray*}
\eqref{sumep} is true since the mean of the left-hand side of
\eqref{sumep} is bounded.

Since $|A(x)|$ is a regularly varying function with index $\rho<0$,
that is
\begin{equation}\label{regvar}
\lim_{t\to\infty}\frac{|A(tx)|}{|A(t)|}=x^{\rho}, ~~x>0.
\end{equation}
It is known that $|A(x)|$ can be written as  $|A(x)|=c(x)f(x)$,
where $\lim_{x\to\infty}c(x)=c>0$ and $f(x)$ is a continuous and
strictly decreasing function on $(0,\infty)$. This implies
\[
\max_{1\le i\le k_n}|A(\frac{m_i}{r_i})|=O(\max_{1\le i\le
k_n}f(\frac{m_i}{r_i}))=O(f(\min_{1\le i\le
k_n}\frac{m_i}{r_i}))=O(f(q_n))=O(|A(q_n)|).
\]
Then it follows from \eqref{RepXPP}, \eqref{sumep} and
Lemma~\ref{lem-sum} that
\begin{eqnarray*}
&&\frac{1}{\sum^{k_n}_{i=1}r_i}\sum^{k_n}_{i=1}r_i\log(X_{m_i,r_i+1}^{(i)})-\log
x_p-\frac{\sum^{k_n}_{i=1}r_ia(m_i,r_i,p)}{\sum^{k_n}_{i=1}r_i}\gamma\\
&=&\frac{\gamma}{\sum^{k_n}_{i=1}r_i}\sum^{k_n}_{i=1}r_i(E_{m_i,r_i+1}^{(i)}-\sum^{m_n}_{j=r_i+1}\frac1j)
   +\frac1{\sum^{k_n}_{i=1}r_i}\sum^{k_n}_{i=1}r_iA(\frac{m_i}{r_i})
\varepsilon(\frac{m_i}{r_i},\frac{e^{E_{m_i,r_i+1}^{(i)}}}{m_i/r_i},\frac{r_i}{m_ip_n})\\
&\le&\frac{\gamma}{\sum^{k_n}_{i=1}r_i}\sum^{k_n}_{i=1}r_i(E_{m_i,r_i+1}^{(i)}-\sum^{m_n}_{j=r_i+1}\frac1j)
   +\frac{O(|A(q_n)|)}{\sum^{k_n}_{i=1}r_i}\sum^{k_n}_{i=1}r_i
|\varepsilon(\frac{m_i}{r_i},\frac{e^{E_{m_i,r_i+1}^{(i)}}}{m_i/r_i},\frac{r_i}{m_ip_n})|\\
&=&O_p(\frac{1}{(\sum^{k_n}_{i=1}r_i)^{1/2}}+|A(q_n)|),
\end{eqnarray*}
proving the lemma. \qed

\vspace{10pt}

\noindent{\it Proof of Theorem~\ref{thm4}}. It follows from
Lemma~\ref{lem3} that
\begin{eqnarray*}
\log\widehat{x}^*_p-\log
x_p&=&\frac{1}{\sum^{k_n}_{j=1}r_j}\sum^{k_n}_{i=1}r_i\log(X_{m_i,r_i+1}^{(i)})-\log
x_p-a_n(p_n)\widehat{\gamma}^*_n\\
&=&-a_n(p_n)(\widehat{\gamma}^*_n-\gamma)+O_p(\frac{1}{(\sum^{k_n}_{i=1}r_i)^{1/2}}+|A(q_n)|),
\end{eqnarray*}
which yields
\begin{equation}\label{share}
\frac{(\sum^{k_n}_{i=1}r_i)^{1/2}}{-a_n(p_n)}\big(\log\widehat{x}^*_p-\log
x_p\big)=
(\sum^{k_n}_{i=1}r_i)^{1/2}(\widehat{\gamma}^*_n-\gamma)+O_p(\frac{1}{-a_n(p_n)}+\frac{(\sum^{k_n}_{i=1}r_i)^{1/2}}{-a_n(p_n)}|A(q_n)|).
\end{equation}
The big ``O" term above converges to zero in probability since
$-a_n(p_n)\to \infty$ from Lemma~\ref{lem2} and
$(\sum^{k_n}_{i=1}r_i)^{1/2}|A(q_n)|\to 0$ as a given condition in
Theorem~\ref{thm3}.  Therefore, the left-hand side of the above
equation converges in distribution to $N(0,\gamma^2)$ by using
Theorem~\ref{thm3}. This completes the proof of Theorem~\ref{thm4}.
\qed

\vspace{10pt}

\noindent{\it Proof of Theorem~\ref{thm2}}. Theorem~\ref{thm2} is
the special case of Theorem~\ref{thm4} except we allow a non-zero
bias term in the limiting distribution. Under the setup in
Theorem~\ref{thm2}, we have $m_i=m_n\sim \frac{n}{k_n}$, and $r_i=r$
is a fixed integer. In the proof of Theorem~\ref{thm4} we have
obtained \eqref{share}. We note that the left-hand side of
\eqref{share} is equal to the left-hand side of \eqref{qnorm}, and
$|a_n(p_n)|=|a(m_n, r, p_n)|\to\infty$. Theorem~\ref{thm1} together
with \eqref{share} yields Theorem~\ref{thm2} if we can show that
$\sqrt{k_n}|A(\frac{m_n}{r})|$ has a finite limit. In fact, we have
from \eqref{regvar} that
\begin{equation}\label{ratio}
\frac{k_n^{1/2}|A(m_n/r)|}{k_n^{1/2}|A(m_n)|}=\frac{|A(m_n/r)|}{|A(m_n)|}\to
r^{-\rho},
\end{equation}
 which coupled with assumption
$k_n^{1/2}A(m_n)\to\delta\in(-\infty, \infty)$ implies
$k_n^{1/2}|A(m_n/r)|\to |\delta|r^{-\rho}$.  This completes the
proof of Theorem~\ref{thm2}. \qed

\vspace{20pt}

\noindent{\it Proof of Theorem~\ref{thm5}.}  In this proof, we will
simply use $m$ and $p$ to denote $m_n$ and $p_n$, respectively.

Define
\[
Z_j^{(i)}=j(\log X_{m,j}^{(i)}-\log X_{m,j+1}^{(i)})
\]
for $j=1,\cdots, r$ and $i=1,\cdots,k_n$.  We have
\[
z_j^{(i)}(y)=Z_j^{(i)}-\frac{1}{a(m,r,p)}(\log(X_{m,r+1}^{(i)})-y)
\]
for $j=1,\cdots, r$ and $i=1,\cdots,k_n$.

Note that we have assumed that $y_0$ is the true value of $\log x_p$.  Now we also assume that $\gamma_0$ is the true value of $\gamma$.
 It follow from the proof of Theorem 4 in Qi~\cite{Qi2010} that
\begin{equation}\label{maxZ}
\max_{1\le j\le r}\max_{1\le i\le
k_n}|Z_j^{(i)}-\gamma_0|=o_p(k_n^{1/2})
\end{equation}
and
\begin{equation}\label{sampleZ}
s_n^2:=\frac{1}{rk_n}\sum^{k_n}_{i=1}\sum^r_{j=1}(Z_j^{(i)}-\gamma_0)^2\stackrel{p}{\to}
\gamma^2_0.
\end{equation}

From now on we will write $z_j^{(i)}(y_0)$ as $z_j^{(i)}$ for convenience.
It follows from \eqref{RepXPP} that
\begin{eqnarray}\label{decomp}
z_j^{(i)}&=&Z_j^{(i)}-\gamma_0+\frac{1}{-a(m,r,p)}(\log(X_{m,r+1}^{(i)})-y_0)+\gamma_0\nonumber\\
&=&(Z_j^{(i)}-\gamma_0)+
\frac{\gamma_0}{-a(m,r,p)}(E_{m,r+1}^{(i)}-\sum^m_{j=r+1}\frac1j)\nonumber\\
&& +\frac{A(m/r)}{-a(m,r,p)}\varepsilon(m/r,\frac{e^{E_{m,r+1}^{(i)}}}{m/r},\frac{r}{mp})\nonumber\\
&=&:
a_j^{(i)}+b_i+c_i.
\end{eqnarray}
We need to show the following three equations:
\begin{equation}\label{maxzy}
\max_{1\le i\le k_n,1\le j\le r}|z_j^{(i)}|=o_p(k_n^{1/2}),
\end{equation}
\begin{equation}\label{samplevar}
\frac{1}{rk_n}\sum^{k_n}_{i=1}\sum^r_{j=1}(z_j^{(i)})^2\stackrel{p}{\to}
\gamma^2_0,
\end{equation}
\begin{equation}\label{CLT}
\frac{1}{\sqrt{rk_n}}\sum^{k_n}_{i=1}\sum^r_{j=1}z_j^{(i)}\stackrel{d}{\to}
N(0,\gamma^2_0).
\end{equation}

Using \eqref{maxep}, \eqref{sumep} with $\ell=2$, \eqref{ratio} and
the fact that $-a(m,r,p)\to\infty$ as $n\to\infty$ from
Lemma~\ref{lem2} we have
\begin{equation}\label{fact1}
\frac{1}{k_n^{1/2}}\max_{1\le i\le k_n}|c_i|\tp 0
~~\mbox{and}~~\frac1{k_n}\sum_{1\le i\le k_n}c_i^2\tp 0.
\end{equation}
We can show
\begin{equation}\label{fact2}
\frac{1}{k_n^{1/2}}\max_{1\le i\le k_n}|b_i|\tp 0
~~\mbox{and}~~\frac1{k_n}\sum_{1\le i\le k_n}b_i^2\tp 0.
\end{equation}
The second expression can be proved by using the estimation that
\[
\frac1{k_n}\E(\sum_{1\le i\le
k_n}b_i^2)\le\frac{\gamma_0^2}{k_n(a(m,r,p))^2}\frac{k_n}{r}=\frac{\gamma_0^2}{r(a(m,r,p))^2}\to
0
\]
from \eqref{mean-var}, and the first one follows from the second one
since
 \[
 \frac{1}{k_n^{1/2}}\max_{1\le i\le k_n}|b_i|\le
\big(\frac1{k_n}\sum_{1\le i\le k_n}b_i^2\big)^{1/2}.
\]

We see that \eqref{maxzy} follows from \eqref{maxZ} and the first
expressions in both \eqref{fact1} and \eqref{fact2}. \eqref{CLT}
follows from Theorem~\ref{thm2} with $\delta=0$ since
\[
\frac{1}{\sqrt{k_nr}}\sum^{k_n}_{i=1}\sum^r_{j=1}z_j^{(i)}=\frac{\sqrt{k_nr}}{-a(m,r,p)}(\log\widehat{x}_p-\log
x_p).
\]

Set $d_i=b_i+c_i$. We have from the Cauchy-Schwarz inequality that
\begin{eqnarray*}
\frac{1}{k_n}\sum^{k_n}_{i=1}d_i^2&=&\frac{1}{k_n}\Big(\sum^{k_n}_{i=1}b_i^2+\sum^{k_n}_{i=1}c_i^2+2\sum^{k_n}_{i=1}b_ic_i\Big)\\
&\le&\frac{1}{k_n}\sum^{k_n}_{i=1}b_i^2+\frac{1}{k_n}\sum^{k_n}_{i=1}c_i^2+2\sqrt{\frac{1}{k_n}\sum^{k_n}_{i=1}b_i^2}\sqrt{\frac{1}{k_n}\sum^{k_n}_{i=1}c_i^2}\\
&\stackrel{p}{\to}&0
\end{eqnarray*}
by using the second expressions in \eqref{fact1} and \eqref{fact2}.
 Now we have from \eqref{decomp} that
\begin{eqnarray*}
\frac{1}{rk_n}\sum^{k_n}_{i=1}\sum^r_{j=1}(z_j^{(i)})^2&=&\frac{1}{rk_n}\sum^{k_n}_{i=1}\sum^r_{j=1}(a_j^{(i)})^2+\frac{1}{k_n}\sum^{k_n}_{i=1}d_i^2+
\frac{2}{rk_n}\sum^{k_n}_{i=1}\sum^r_{j=1}a_j^{(i)}d_i.
\end{eqnarray*}
On the right-hand side in the above equation, the first term
converges in probability to $\gamma_0^2$ from \eqref{sampleZ}, the
second term converges in probability to zero, and the third term
converges in probability to zero by using the Cauchy-Schwarz
inequality. This completes the proof of \eqref{samplevar}.

The proof for \eqref{wilks} is quite standard under conditions
\eqref{maxzy}, \eqref{samplevar} and \eqref{CLT}; see e.g.,
Owen~\cite{Owen2001} for details. \qed


\vspace{20pt}

\noindent{\it Proof of Theorem~\ref{thm6}.}  By following the same
arguments in the proof of Theorem~\ref{thm5}, it suffices to verify
the following three conditions
\begin{equation}\label{maxzy1}
\max\Big(\max_{1\le i\le k_n,1\le j\le r}|z_j^{(i)}|, |z(y_0)|\Big)=o_p(k_n^{1/2}),
\end{equation}
\begin{equation}\label{samplevar1}
\frac{1}{rk_n+1}\Big(\sum^{k_n}_{i=1}\sum^r_{j=1}(z_j^{(i)})^2+z(y_0)^2\Big)\stackrel{p}{\to}
\gamma^2_0,
\end{equation}
\begin{equation}\label{CLT1}
\frac{1}{\sqrt{rk_n+1}}\Big(\sum^{k_n}_{i=1}\sum^r_{j=1}z_j^{(i)}+z(y_0)\Big)\stackrel{d}{\to}
N(0,\gamma^2_0).
\end{equation}
\eqref{maxzy1}, \eqref{samplevar1} and \eqref{CLT1} follow from \eqref{maxzy}, \eqref{samplevar} and \eqref{CLT} since
\[
z(y_0)=-\frac{a_n}{k_nr}\sum^{k_n}_{i=1}\sum^r_{j=1}z_j^{(i)}(y_0)=O_p(\frac{a_n}{\sqrt{k_n}})=o_p(k_n^{1/6})
\]
from \eqref{zy}. This completes the proof of Theorem~\ref{thm6}. \qed



\begin{thebibliography}{AA} 





\bibitem{AEG2022} Allouche, M., El Methni, J. and Girard, S. (2022). A refined Weissman
estimator for extreme quantiles. {\em Extremes}.
https://doi.org/10.1007/s10687-022-00452-8


\bibitem{BGT1987}
Bingham, N. H.,  Goldie, C. M. and Teugels, J. L. (1987). {\em
Regular Variation}.  New York: Cambridge University.


\bibitem{CVA2008}
Chen, J., Variyath, A. M.  and Abraham, B. (2008). Adjusted
empirical likelihood and its properties. {\em Journal of
Computational and Graphical Statistics}, {\bf 17}, 426--443.





\bibitem{DPR2000}
Davydov, Yu., Paulauskas, V. and Ra\v ckauskas, A.(2000). More on
$P$-stable convex sets in Banach spaces. {\em Journal of Theoretical
Probability}, {\bf 13}, 39--64.

\bibitem{deHannF2006}
De Haan, L. and Ferreira, A. (2006). Extreme Value Theory: An Introduction.  Springer, New York.


\bibitem{DR1993}
De Haan, L. and Rootz\'n, H. (1993). On the estimation of high
quantiles. {\em Journal of Statistical Planning and Inference}, {\bf
35}, 1--13.

\bibitem{deHannS1996}
De Haan, L. and Stadtm\"{u}ller, U.  (1996). Generalized regular
variation of second order. {\em Journal of the Australian
Mathematical Society (Series A)}, {\bf 61}, 381--395.





\bibitem{DH1989}
Dekkers, A. L. M. and de Haan, L. (1989). On the estimation of the
extreme-value index and large quantile estimation. {\em  The Annals
of Statistics}, {\bf 17}, 1795--1832.

\bibitem{DEH1989}
Dekkers, A. L. M., Einmahl, J. H. J. and de Haan, L. (1989). A
moment estimator for the index of an extreme-value distribution.
{\em The Annals of Statistics}, {\bf 17}, 1833--1855.



\bibitem{ELM1997}
Embrechts, P., Kl\"uppelberg, C. and Mikosch, T. (1997). {\em
Modelling Extremal Events for Insurance and Finance}. Berlin:
Springer.


\bibitem{FDP2003}
Ferreira, A., de Haan, L. and  Peng, L. (2003). On optimising the
estimation of high quantiles of a probability distribution. {\em
Statistics}, {\bf 37}, 401--434.




\bibitem{GHP2002}
Gomes, M.I., de Haan, L. and Peng, L. (2002). Semi-parametric
estimation of the second order parameter in statistics of extremes.
{\em Extremes} 5, 387--414.


\bibitem{GLPY2015}
 Gong, J.,  Li, Y., Peng, L. and Yao, Q. (2015).
Estimation of extreme quantiles for functions of dependent random
variables. {\em Journal of the Royal Statistical Society B}, {\bf
77}, 1001--1024.



\bibitem{HL1990}
Hall, P.  and La Scala, L. (1990). Methodology and algorithms of
empirical likelihood. {\em International Statistical Review}, {\bf
58}, 109--127.


\bibitem{Hill1975}
Hill, B. M. (1975). A simple general approach to inference about the
tail of a distribution. {\em The Annals of Statistics}, {\bf 3},
1163--1174.

\bibitem{HPS2022}
Hu, S., Peng Z., and Nadarajah, S. (2022). Location invariant heavy
tail index estimation with block method.  {\em Statistics}, {\bf
56},  479-497.




\bibitem{LiQi2019}
Li, Y. and Qi, Y. (2019). Adjusted empirical likelihood method for
the tail index of a heavy-tailed distribution. {\em Statistics and
Probability Letters}, {\bf 152}, 50--58.


\bibitem{LC2010}
Liu, Y. and Chen, J. (2010). Adjusted empirical likelihood with
high-order precision. {\em Annals of Statistics}, {\bf 38},
1341--1362.


\bibitem{LP2002}
Lu, J. and Peng, L. (2002). Likelihood based confidence intervals
for the tail index. {\em Extremes}, {\bf 5(4)}, 337--352



\bibitem{Owen1988}
Owen, A. (1988). Empirical likelihood ratio confidence intervals for
a single functional. {\em Biometrika},  {\bf 75}, 237--249.

\bibitem{Owen1990}
Owen, A. (1990). Empirical likelihood regions. {\em The Annals of
Statistics}, {\bf 18}, 90--120.


\bibitem{Owen2001}
Owen, A. (2001). {\em Empirical Likelihood.} Chapman and Hall.


\bibitem{Paulauska2003}
Paulauskas, V. (2003). A new estimator for a tail index, {\em Acta
Applicandae Mathematicae}, {\bf 79}, 55--67.

\bibitem{PV2011}
Paulauskas, V. and Vaiciulis, M. (2011). Several modifications of
DPR estimator of the tail index. {\em Lithuanian Mathematical
Journal}, {\bf 51}, 36--50.

\bibitem{PV2012}
Paulauskas, V. and Vaiciulis, M. (2012). Estimation of the tail
index in the max-aggregation scheme. {\em Lithuanian Mathematical
Journal}, {\bf 52}, 297--315.



\bibitem{PV2013}
Paulauskas, V. and Vaiciulis, M. (2013). On an improvement of Hill
and some other estimators. {\em Lithuanian Mathematical Journal},
{\bf 53}, 336--355.

\bibitem{PV2017}
Paulauskas, V. and Vaiciulis, M. (2017). A class of new tail index
estimators.  {\em Annals of the Institute of Statistical
Mathematics}, \textbf{69}, 461--487.



\bibitem{PQ2004}
Peng, L. and  Qi, Y. (2004).  Estimating the first and second order
parameters of a heavy tailed distribution.  {\em Australian \& New
Zealand Journal of Statistics} 46, no 2, 305--312.


\bibitem{PengQi2006a}
Peng, L. and Qi, Y. (2006a). A new calibration method of
 constructing empirical likelihood-based confidence intervals for the tail
 index. {\em Australian and New Zealand Journal of Statistics}, {\bf 48},
 59--66.

\bibitem{PengQi2006b}
Peng, L. and Qi, Y. (2006b). Confidence regions for high quantiles
of a heavy-tailed distribution.  {\em The Annals of Statistics},
{\bf 34} (4), 1964--1986.

\bibitem{PengQi2017}
Peng, L. and Qi, Y. (2017). {\em Inference for Heavy-Tailed Data: Applications in Insurance and Finance}
Aacdemic Press.



\bibitem{Pickands1975}
Pickands, J. (1975). Statistical inference using extreme order
statistics. {\em The Annals of Statistics}, {\bf 3}, 119--131.


\bibitem{Qi2010}
Qi, Y. (2010). On the tail index of a heavy tailed distribution.
{\em Annals of the Institute of Statistical Mathematics}, {\bf 62},
277--298.




\bibitem{QL1994}
Qin, J. and Lawless, J. (1994). Empirical likelihood and general
estimating equations. {\em The Annals of Statistics}, {\bf 22},
300--325.



\bibitem{Tsao2004}
Tsao, M. (2004). A new method of calibration for the empirical
loglikelihood ratio. {\em Statistics and Probability Letters}, {\bf
68}, 305--314.

\bibitem{Vai2012}
Vaiciulis, M. (2012). Asymptotic properties of generalized DPR
statistic. {\em Lithuanian Mathematical Journal}, {\bf 52}, 95--110.

\bibitem{Vai2014}
Vaiciulis, M. (2014). Local-maximum-based tail index estimator. {\em
Lithuanian Mathematical Journal}, {\bf 54}, 503--526.





\bibitem{Weissman1978}
Weissman, I. (1978). Estimation of parameters and large quantiles
based on the k largest observations. {\em J. Am. Stat. Assoc.}, {\bf
73}, 812--815.

\bibitem{XP2020}
Xiong, L. and Peng, Z. (2020). Heavy tail index estimation based on
block order statistics.  {\em Journal of Statistical Computation and
Simulation}, {\bf 90}, 2198--2208.


\end{thebibliography}

\vspace{20pt}

\noindent\textbf{Acknowledgements}.

The authors would like to thank two anonymous referees for their
constructive suggestions that led to improvement of the paper. The
research of Yongcheng Qi was supported in part by NSF Grant
DMS-1916014. The research of Jingping Yang was supported by the
National Natural Science Foundation of China (Grants No. 12071016).

\noindent\textbf{Disclosure statement}.

The authors report there are no competing interests to declare.

\baselineskip 12pt
\def\ref{\par\noindent\hangindent 25pt}

\end{document}